\documentclass[11pt,thmsa,emstex]{article}
\usepackage[left=3cm,right=3cm,top=3cm,bottom=3cm]{geometry}
\usepackage{amsfonts}
\usepackage{t1enc}
\usepackage{amssymb}
\usepackage{epsfig}
\usepackage[english]{babel}
\usepackage{amsmath}
\usepackage{graphics}
\usepackage{layout}
\usepackage{latexsym}
\usepackage{tikz}
\usepackage{amsmath,amsxtra,amssymb,latexsym,epsfig,amscd,amsthm,fancybox,epsfig}
\usepackage{authblk}

\usepackage{hyperref}

\newtheorem{theorem}{Theorem}[section]
\newtheorem{lemma}[theorem]{Lemma}
\newtheorem{corollary}[theorem]{Corollary}
\newtheorem{proposition}[theorem]{Proposition}
\newtheorem{example}[theorem]{Example}
\newtheorem{remark}[theorem]{Remark}

\newtheorem{definition}[theorem]{Definition}
\newtheorem{exercice}[theorem]{Exercice}

\def\bit{\begin{itemize}}

\def\eit{\end{itemize}}

\reversemarginpar   

\def\bc{\begin{center}}

\def\ec{\end{center}}

\def\bthm{\begin{theorem}}

\def\ethm{\end{theorem}}

\def\bcor{\begin{corollary}}

\def\ecor{\end{corollary}}

\def\bprop{\begin{proposition}}

\def\eprop{\end{proposition}}

\def\blem{\begin{lemma}}

\def\elem{\end{lemma}}

\def\bex{\begin{example}}

\def\eex{\end{example}}

\def\bexo{\begin{exercice} \rm }

\def\eexo{\end{exercice} }

\def\brem{\begin{remark}}

\def\erem{\end{remark}}

\def\prf{{\bf Proof }}

\def\bdes{\begin{description}}

\def\edes{\end{description}}

\def\beq{\begin{equation}}

\def\eeq{\end{equation}}

\def\ben{\begin{enumerate}}

\def\een{\end{enumerate}}

\def\beqar{\begin{eqnarray}}

\def\eeqar{\end{eqnarray}}

\def\beqarr{\begin{eqnarray*}}

\def\eeqarr{\end{eqnarray*}}

\def\qed{\hspace{.1in}{\bf QED}\\[2ex]}

\def\prf{{\bf Proof }\hspace{.1in}}

\def\Pr{{\mathsf P}}

\def\ZZ{{\mathbb Z}}       
\def\RR{{\mathbb R}}  

\def\NN{{\mathbb N}}

\def\TT{{\mathbb T}}
\def\ZZ{{\mathbb Z}}

\def\eps{\varepsilon}

\newcommand{\bi}{\mathbf{i}}
\newcommand{\bu}{\mathbf{u}}
\newcommand{\bv}{\mathbf{v}}
\newcommand{\bphi}{\mathbf{\Phi}}

\def\1{{\rm 1\mskip-4.4mu l}}

\usetikzlibrary{decorations.markings}
\tikzset{degil/.style={
            decoration={markings,
            mark= at position 0.5 with {
                  \node[transform shape] (tempnode) {$\backslash$};
                  }
              },
              postaction={decorate}
}
}
\begin{document}
\title{A user-friendly condition for exponential ergodicity in randomly switched environments}
\author[1]{Michel Bena{\"i}m }

\author[2]{Tobias Hurth }
\author[1]{Edouard Strickler}
\affil[1]{ Institut de Math\'ematiques\\Universit\'e de Neuch\^atel, Switzerland}
\affil[2]{ Ecole Polytechnique F\'ed\'erale de Lausanne, Switzerland}
\maketitle
\begin{abstract}
We consider random switching between finitely many vector fields leaving positively invariant a compact set. Recently, Li, Liu and Cui showed in \cite{li17} that if one the vector fields has a globally asymptotically stable (G.A.S.) equilibrium from which one can reach a point satisfying a weak H{\"o}rmander-bracket condition, then the process converges in total variation to a unique invariant probability measure. In this note, adapting the proof in \cite{li17} and using results of \cite{BMZIHP}, the assumption of a G.A.S. equilibrium is weakened to the existence of an accessible point at which a barycentric combination of the vector fields vanishes. Some  examples are given which demonstrate the usefulness of this condition.
\end{abstract}


\section{Introduction}

Let $E=\{1,\ldots,N\}$ be a finite set and  $\mathrm{F}  = \{F^i\}_{i \in E}$ a family of smooth globally integrable vector fields on $\RR^d.$
For each $i \in E$ we let $\varphi^i = \{\varphi^i_t\}$   denote the  flow induced by $F^i.$
We assume  throughout that there exists a  compact set $M \subset \RR^d$  which is {\em positively invariant} under each $\varphi^i.$ That is $$\varphi^i_t(M) \subset M$$ for all $t \geq 0$.  Our assumption that $M \subset \RR^d$ is mostly for convenience. The results of this note can readily be generalized to the situation where $M$ is a subset of a finite-dimensional smooth manifold.  

Consider a Markov process $Z = (Z_t)_{t \geq 0}, Z_t = (X_t, I_t),$  living on $M \times E$ whose infinitesimal generator acts on functions $g : M \times E \mapsto \RR,$ smooth in the first variable, according to the formula
\beq
\label{eq:defcL}
{\cal L}g(x,i) =  \langle F^i(x), \nabla g^i(x) \rangle + \sum_{j \in E} a_{ij}(x) (g^j(x) - g^i(x)),
\eeq
where $g^i(x)$ stands for $g(x,i)$ and  $a(x) = (a_{ij}(x))_{i,j \in E}$ is an irreducible {\em rate} matrix continuous in $x.$ Here, by  rate matrix, we mean a matrix  having nonnegative off diagonal entries and zero diagonal entries.

In other words, the dynamics of $X$ is given by an ordinary differential equation
\beq
 \frac{dX_t}{dt} = F^{I_t}(X_t),
\label{eq:pdmp}
\eeq
while $I$ is a continuous-time jump process taking values in $E$ controlled by $X:$
$$\Pr(I_{t+s} = j | {\cal F}_t)  = a_{ij}(X_t) s + o(s) \mbox{ for } j \neq i \mbox{ on } \{I_t = i\},$$
where ${\cal F}_t = \sigma ((X_s,I_s) \: : s \leq t\}.$

The process $Z$ belongs to the class of processes called Piecewise Deterministic Markov Processes (PDMP), introduced by Davis in \cite{Dav84}. Ergodic properties of these processes have recently been the focus of much attention (e.g., \cite{Le_Borgne}, \cite{Cloez},\cite{BMZIHP}, \cite{BCL17}, \cite{BS17}, \cite{torus}).

Using the terminology in \cite{BMZIHP}, a point $x^{\star} \in M$ is said to satisfy the \textit{weak bracket condition} if the Lie algebra generated by  $(F^i)_{i \in E}$ at $x^{\star}$ has full rank. If such a point is furthermore accessible (meaning that every neighborhood of $x^{\star}$ is reached with positive probability by $X_t$), then the process admits a unique invariant probability measure which is absolutely continuous with respect to the Lebesgue measure on $M \times E$ (see e.g \cite[Theorem 1]{BH12} or \cite[Theorem 4.5]{BMZIHP}). If the weak bracket condition is replaced by the so-called \textit{strong bracket condition} (cf. Definition \ref{def:bracket_condition} below), the process then converges in total variation (see \cite[Theorem 4.6]{BMZIHP}). Simple examples show that the weak bracket condition itself is not sufficient to ensure convergence (cf. \cite{BH12}).

Recently, Li, Liu and Cui showed in \cite{li17} that the two following conditions yield convergence in total variation (see \cite[Theorem 9]{li17}) : 
\begin{enumerate}
\item[\textbf{(i')}] There exists a  globally asymptotically stable (G.A.S.) equilibrium for one of the vector fields,
\item[\textbf{(ii)}] The weak bracket condition holds at an accessible point.
\end{enumerate}
    In this note we replace  \textbf{(i')} by the more general condition
    \begin{enumerate}
\item[\textbf{(i)}] There exists  an accessible point $e^{\star}$ at which a barycentric combination of the vector fields vanishes,
\end{enumerate}
and  prove exponential convergence in total variation (see Theorem \ref{main} and Corollary \ref{cor:cvVT}).  Our proof is inspired by \cite{li17}  but is simplified using results of \cite{BMZIHP}. 

It turns out that when the vector fields are analytic, \textbf{(i)} and \textbf{(ii)} imply the strong bracket condition at  $e^{\star}$ (cf. Proposition \ref{prop:strong_bracket_analytic}). Nonetheless,  \textbf{(i)} and  \textbf{(ii)} are usually much easier to verify than the strong bracket condition. This is  illustrated by the examples in Section \ref{sec:app}. In the nonanalytic case, neither condition implies the other as shown in Section \ref{sec:strong_condition} (see Examples \ref{ex:nonanalytic} and \ref{ex:nonzero}). All these results are summarized in the following scheme.

\bigskip

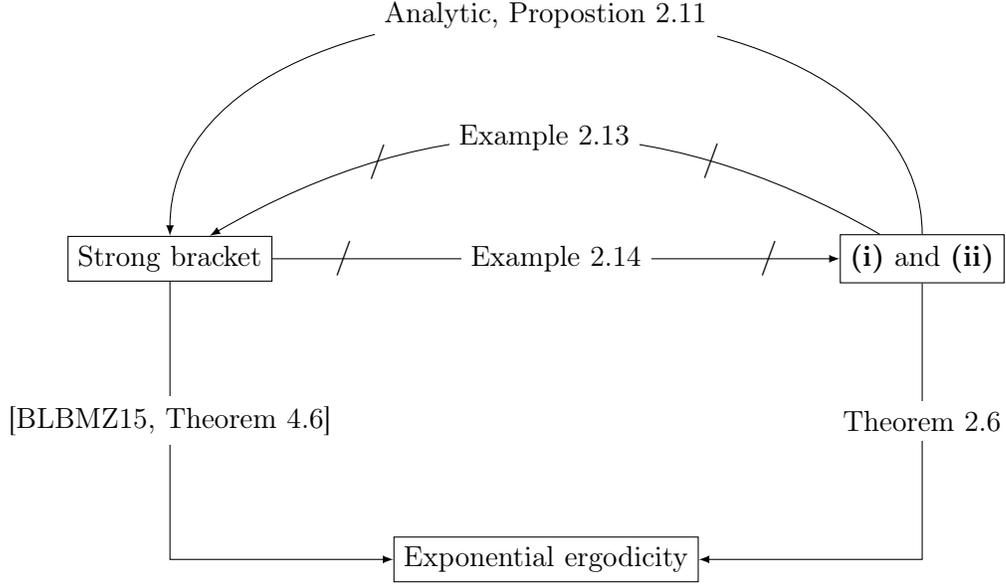
\begin{tikzpicture}
\node[draw] (P) at (0,0) {Strong bracket};
\node[draw] (S) at (10,0) {\textbf{(i)} and \textbf{(ii)}};
\node[draw] (N) at (5,-4) {Exponential ergodicity};
\draw[->,>=latex] (P) -- node[very near start]{\big{/}} node[midway,fill=white]{Example \ref{ex:nonzero}} node[very near end]{\big{/}} (S);
\draw[->,>=latex] (P) |- node[near start,fill=white]{\cite[Theorem 4.6]{BMZIHP}} (N);
\draw[->,>=latex] (S) |- node[near start,fill=white]{Theorem \ref{main}} (N) ;
\draw[->,>=latex] (S) to [bend right=90]  node[midway,fill=white]{Analytic, Propostion \ref{prop:strong_bracket_analytic} } (P);
\draw[->,>=latex] (S) to [bend right=30]  node[near start]{\big{/}} node[midway,fill=white]{Example \ref{ex:nonanalytic}}   node[near end]{\big{/}} (P);
\end{tikzpicture}

\bigskip


\section{Definitions and main results}
\label{sec:def}

We begin by recalling some general definitions. Let $(P_t)_{t \geq 0}$ be a Markov semigroup  on a metric space $\mathcal{M}$.  
\begin{definition}
We say that $z^* \in \mathcal{M}$ is a \emph{Doeblin point} if there exists a neighborhood $U$ of $z^*$, a nonzero measure $\nu$ and positive real numbers $t^*,c$ such that $P_{t^*}(z,\cdot) \geq c \nu( \cdot)$ for all $z \in U$. 
\end{definition}

\begin{definition}      \label{def:accessible}
We say that $z^* \in \mathcal{M}$ is $(P_t)$-\emph{accessible} from $B \subset \mathcal{M}$ if for every neighborhood $U$ of $z^*$ and for all $z \in B$, there exists a positive real $t$ such that $P_{t}(z,U) > 0$. 
\end{definition}

In the specific context of PDMPs, the latter definition can be expressed more intuitively as follows.  For  $\bi=(i_1, \ldots, i_m) \in E^m$ and $\bu = (u_1,\ldots,u_m) \in \RR_+^m$, we denote by $\bphi_{\bu}^{\bi}$ the composite flow : $\bphi_{\bu}^{\bi} = \varphi_{u_m}^{i_m} \circ \ldots \circ \varphi_{u_1}^{i_1}$.  For $x \in M$ and $t \geq 0$, we denote by $\gamma^+_t(x)$ (resp. $\gamma^+(x)$)  the set of points that are reachable from $x$ at time $t$ (resp. at any nonnegative time) with a composite flow: $$\gamma^+_t(x)=\{ \bphi_{\bv}^{\bi}(x), \: (\bi,\bv) \in E^m \times \RR_+^m, m \in \NN, v_1 + \ldots + v_m = t\},$$

$$ 
\gamma^+(x) = \bigcup_{t \geq 0} \gamma^+_t(x).
$$ 

\begin{definition}       \label{def:F_accessible}
A point $x^* \in M$ is $\{F^i\}$-\emph{accessible} from $B \subset M$ if $x^* \in \cap_{x \in B} \overline{\gamma^+(x)}$.
\end{definition}   
From now on, we let $(P_t)_{t \geq 0}$ be the semigroup induced by $(Z_t)_{t \geq 0}$ on $\mathcal{M}= M \times E$. 
Because of the irreducibility assumption on the rate matrix $a(x)$, Definitions~\ref{def:accessible} and~\ref{def:F_accessible} co\"incide (see e.g. \cite[Lemma 3.2]{BMZIHP},  or \cite[Lemma 3.1]{BCL17}):
\bprop
For all $j, k \in E$, the point $(x^*,j) \in M \times E$ is $(P_t)$-accessible from $B \times \{k\} \subset M \times E$ if and only if $x^*$ is $\{F^i\}$-accessible from $B$. 
\label{acces}
\eprop 
Therefore, in the sequel, we will say that a point $x^* \in M$ is accessible from $B \subset M$ if it is $\{F^i\}$-accessible from $B$. We will simply say that $x^*$ is accessible if it is $\{F^i\}$-accessible from $M$.  
Set $\mathrm{F}_0 = \{F^i\}_{i \in E}$ , $\mathrm{F}_{k+1} = \mathrm{F}_k \cup \{[F^i, V], V \in \mathrm{F}_k\}$, $\mathcal{F}_{0} = \{F^i  - F^j\: : i, j = 1, \ldots m\}$ and $\mathcal{F}_{k+1}  =  \mathcal{F}_k \cup \{[F^i , V ] \: : V \in \mathcal{F}_k\}.$ Here $[\cdot , \cdot]$ stands for  the Lie bracket operation, which is defined as 
$$ 
[V,W](x) = DW(x) V(x) - DV(x) W(x), \quad x \in \RR^d, 
$$ 
for smooth vector fields $V$ and $W$ on $\RR^d$ with differentials $DV$ and $DW$. The following definition is given in \cite{BMZIHP}.

\begin{definition}
\label{def:bracket_condition}
We say  that  the {\em weak bracket} (resp. {\em strong bracket})  condition holds at $p \in M$ if the vector space spanned by the vectors $\{V(p) \: : V \in \cup_{k \geq 0}  \mathrm{F}_{k}\}$ (resp. $\{V(p) \: : V \in \cup_{k \geq 0}  \mathcal{F}_{k}\}$)  has full rank. 
\end{definition}
It is clear from this definition that the strong bracket condition implies the weak one. Weak bracket and strong bracket conditions are equivalent to Condition~B and Condition~A in~\cite{BH12}, respectively.  The weak bracket condition is closely related to the classical H\"ormander hypoellipticity condition that yields smoothness of transition densities for diffusions (see e.g.~\cite{Nualart}). More background on the weak and strong bracket conditions with an emphasis on how they relate to controllability is provided in~\cite{SJ72}.

\subsection{Main result} 
We now state our main result.
 
\bthm
Suppose that
\begin{enumerate}
\item[\emph{\textbf{(i)}}] There exist $\alpha_1,\ldots,\alpha_N \in \RR$ with $\sum \alpha_i = 1$ and $e^{\star} \in M$ such that  $\sum \alpha_i F^i(e^{\star}) = 0$,
\item[\emph{\textbf{(ii)}}] There exists a point $x^*$ accessible from $\{e^{\star}\}$ where the weak bracket condition holds. \end{enumerate}
Then for all $j \in E$, $(e^{\star},j)$ is a Doeblin point.
\label{main} 
\ethm

Note that we do not impose that the $\alpha_i$ are nonnegative. In particular, condition \textbf{(i)} holds whenever two vector fields at some point are collinear but not equal.

The following corollary is a consequence of standard results (see e.g \cite[Theorem 4.6]{BMZIHP} for a proof).

\bcor \label{cor:cvVT}
In addition to the assumptions in Theorem \ref{main}, suppose that $e^{\star}$ is accessible. Then, the process $Z$ admits a unique invariant probability measure $\pi$ which is absolutely continuous with respect to Lebesgue measure. Moreover, there exist positive constants $C, \gamma$ such that for all $t \geq 0$ and for all $(x,i) \in M \times E$, $$ \| P_t((x,i),\cdot) - \pi \|_{TV} \leq C e^{- \gamma t}.$$
\ecor

In Section \ref{sec:app}, we give more applications in a stochastic persistence context, relying on recent results in \cite{B18}.
Theorem~\ref{main} is a direct consequence of Theorem~4.2 in~\cite{BMZIHP} and of Proposition~\ref{prop:submersion_at_E} that we state below.  For convenience, we also record a version of Theorem~4.2 from~\cite{BMZIHP}. Here and throughout, for $s > 0$ and $m \in \NN^*$, we set $D_m^s = \{ \bv \in \RR_+^m \: : v_1 + \ldots + v_m \leq s \}$. 

\bthm[Bena\"im -- Le Borgne -- Malrieu -- Zitt]
Let $x$ be a point of $M$, $(\bi,\bu)$ and $s > u_1 + \ldots + u_m$ such that the map $\Psi^s : D^s_m \to  \RR^d$, $ \bv \to \varphi^{i_{m+1}}_{s-(v_1+ \ldots + v_m)} \circ \bphi_{\bv}^{\bi}(x)$ is a submersion at $\bu$.  Then for all $j \in E$, $(x,j)$ is a Doeblin point.
\label{thmsubdoeb}
\ethm 

\bprop      \label{prop:submersion_at_E}
Under conditions \textbf{(i)} and \textbf{(ii)} of Theorem \ref{main},  there exist $s > 0$, $i_{m+1} \in E$, $\bi \in E^m$ and $\bu \in \RR^m_+$ with $u_1 + \ldots + u_m < s$ such that the map $\Psi: D^s_m \to  \RR^d$, $ \bv \to \varphi^{i_{m+1}}_{s-(v_1+\ldots+v_m)} \circ \bphi_{\bv}^{\bi}(e^{\star})$ is a submersion at $\bu$. 
\eprop 

\subsection{Links with the strong bracket condition}
\label{sec:strong_condition}
In \cite{BMZIHP} and \cite{BH12}, the authors show that the conclusions of Theorems \ref{main} and \ref{cor:cvVT} hold when the weak bracket condition is replaced by the strong one. A natural question is whether our assumptions already imply that the strong bracket condition holds at some point.  We address this question in Propositions \ref{prop:basic} and \ref{prop:strong_bracket_analytic}.   

\begin{proposition}
\label{prop:basic}
Let $e^{\star} \in M$ satisfy condition \emph{\textbf{(i)}} of Theorem \ref{main}. Suppose further that the weak bracket condition holds at  $e^{\star}$. Then, the strong bracket condition is also satisfied at $e^{\star}$.
\end{proposition}

\prf To simplify notation, we set 
$$
W(e^{\star}) = \{V(e^{\star}): V \in \cup_{k \geq 0} \mathrm{F}_k\}, \quad 
S(e^{\star}) = \{V(e^{\star}): V \in \cup_{k \geq 0} \mathcal{F}_k\}. 
$$
We will show that the linear spans of $W(e^{\star})$ and $S(e^{\star})$ are equal to each other, which then implies the proposition. It is clear that the span of $S(e^{\star})$ is a subspace of the span of $W(e^{\star})$. Therefore, it suffices to show that $W(e^{\star})$ is contained in the span of $S(e^{\star})$.  Fix a vector field $V \in \cup_{k \geq 0} \mathrm{F}_k$ and let $j$ be the smallest nonnegative integer such that $V \in \mathrm{F}_j$. By induction it is not hard to see that for any $i \geq 1$, the collection of vector fields $\mathrm{F}_i \setminus \mathrm{F}_{i-1}$ is contained in the span of $\cup_{k \geq 0} \mathcal{F}_k$. Thus, if $j \geq 1$, the point $V(e^{\star})$ lies in the span of $S(e^{\star})$. If $j=0$, there is $l \in E$ such that $V = F^l$. By condition \textbf{(i)}, there are real numbers $(\alpha_i)_{i \in E}$ such that $\sum_{i \in E} \alpha_i = 1$ and $\sum_{i \in E} \alpha_i F^i(e^{\star}) = 0$.    Therefore, 
$$ 
F^l(e^{\star}) = \sum_{i \in E} \alpha_i  F^l(e^{\star}) - \sum_{i \in E} \alpha_ i F^i(e^{\star}) = \sum_{i \in E} \alpha_i (F^l(e^{\star}) - F^i(e^{\star})).  
$$ 
Since the vector fields $(F^l - F^i)_{i \in E}$ lie in $\mathcal{F}_0$, we have again that $V(e^{\star})$ is in the span of $S(e^{\star})$. This finishes the proof. 
\qed

\bprop                 \label{prop:strong_bracket_analytic}
Assume that for all $i \in E$, $F^i$ is  analytic and that the assumptions of Theorem \ref{main} hold. Then $e^{\star}$ satisfies the strong bracket condition.
\eprop

In most applications, the vector fields governing the PDMP are analytic (see also Section~\ref{sec:app}). As a consequence, the interest of Theorem \ref{main} lies essentially in the fact that the weak bracket condition is easier to verify than the strong one.  The proof of Proposition \ref{prop:strong_bracket_analytic} relies on the following result, due to Sussmann and Jurdjevic ~\cite[Corollary 4.7]{SJ72}.

%

\bthm[Sussmann -- Jurdjevic]   \label{thm:str_acc_bracket}
Assume that the vector fields $(F^i)_{i \in E}$ are analytic, and let $x$ be any point in $M$. Then, there is $t > 0$ such that $\gamma^+_t(x)$ has nonempty interior if and only if the strong bracket condition holds at $x$. 
\ethm

\prf \textbf{of Proposition~\ref{prop:strong_bracket_analytic}} 

By Proposition~\ref{prop:submersion_at_E}, there are $s > 0$, $i_{m+1} \in E$, $\bi \in E^m$ and $\bu \in \RR^m_+$ with $u_1 + \ldots + u_m < s$ such that $\Psi: \bv \to \varphi^{i_{m+1}}_{s-(v_1+\ldots+v_m)} \circ \bphi_{\bv}^{\bi}(e^{\star})$ is a submersion at $\bu$.  By the constant-rank theorem, there exists an open neighborhood $U$ of $\bu$ such that $\Psi(U)$ is open. Without loss of generality, we can assume that $v_1+\ldots+v_m < s$ for all $\bv \in U$. Then, $\Psi(U)$ is a nonempty open subset of $\gamma^+_s(e^{\star})$. By Theorem~\ref{thm:str_acc_bracket}, $e^{\star}$ satisfies the strong bracket condition.  
\qed

From a more theoretical point of view, we now provide an example in the plane where conditions \textbf{(i)} and \textbf{(ii)} are satisfied, but, in the absence of analyticity, there is no point where the strong bracket condition holds.     
 
\bex \rm \label{ex:nonanalytic} We work in polar coordinates $(\theta, r)$. On the annulus 
$$ 
M = \left\{(\theta,r): \tfrac{1}{2} \leq r \leq 2\right\}, 
$$
we switch between vector fields $F^0(\theta,r) = (1,h(r))^T$ and $F^1 (\theta,r) = (f(\theta),g(\theta) + h(r))^T$, where 
$$ 
h(r) = r (1-r),
$$ 
and where $f$ and $g$ satisfy the following properties: 
\begin{enumerate}
\item The functions $f$ and $g$ are $C^{\infty}$ and $2 \pi$-periodic on $\RR$. 
\item We have $0 < f \leq 1$ and $0 \leq g \leq 1$. 
\item We have $f(\tfrac{\pi}{2}) = \tfrac{1}{2}$ and $g(0) > 0$. Moreover, there is  $\epsilon \in (0, \tfrac{\pi}{4})$ such that $f(\theta) = 1$ for $\lvert \theta - \tfrac{\pi}{2} \rvert > \epsilon$ and $g(\theta) = 0$ for $\lvert \theta \rvert > \epsilon$. 
\end{enumerate}
It is easy to see that such functions $f$ and $g$ exist and that they cannot be analytic. Also note that $M$ is positively invariant under the flows associated with $F^0$ and $F^1$ because $h(\tfrac{1}{2}) > 0$ and $g(\theta) + h(2) < 0$ for all $\theta$. Since $M$ is compact and since $f$, $g$ and $h$ are smooth functions, the vector fields $F^0$ and $F^1$ are globally integrable.  

The point $e^{\star} = (\tfrac{\pi}{2},1)^T$ is an equilibrium point of the vector field $2 F^1-F^0$, so condition \textbf{(i)} is satisfied. Since $h(r) > 0$ for $r \in (0,1)$ and $h(r) < 0$ for $r > 1$, the unit circle is a global attractor of $F^0$.  Thus, any point on the unit circle, in particular the point $e^{\star}$, is accessible from any starting point in $M$. The weak bracket condition holds at the point $(0,1)^T$ because $F^0(0,1) = (1,0)^T$ and $F^1(0,1) = (1,g(0))^T$ generate the entire tangent space at $(0,1)^T$. As $(0,1)^T$ lies on the unit circle, it is accessible from $e^{\star}$.  

It remains to show that the strong bracket condition is nowhere satisfied.  We have 
$$ 
[F^0, F^1](\theta,r) = (f'(\theta), g'(\theta) - h'(r) g(\theta))^T 
$$ 
and 
$$ 
F^1(\theta,r) - F^0(\theta,r) = (f(\theta)-1, g(\theta))^T. 
$$ 
If $\lvert \theta - \tfrac{\pi}{2} \rvert > \epsilon$, both $[F^0,F^1](\theta,r)$ and $(F^1-F^0)(\theta,r)$ have $\theta$-coordinate $0$. And if $\lvert \theta \rvert > \epsilon$, the $r$-coordinate of $[F^0,F^1]$ and $F^1 - F^0$ vanishes.  Now, let $k(\theta,r)$ be a smooth function and let $K_i(\theta,r) = k(\theta,r) (1-i,i)^T$ for $i \in \{0,1\}$.  Then, 
\begin{align*}
[F^0,K_1](\theta,r) =& (0,\ast)^T, & [F^0,K_0](\theta,r) =& (\ast,0)^T, \\
 [F^1,K_1](\theta,r) =& (0,\ast)^T, & [F^1,K_0](\theta,r) =& (\ast, -g'(\theta) k(\theta,r))^T,
\end{align*}
and $g'(\theta) k(\theta,r) = 0$ for $\lvert \theta \rvert > \epsilon$. Here, $\ast$ stands for some term, possibly depending on $\theta$ and $r$, that may differ from equation to equation. This shows that for any $(\theta,r) \in M$, $V(\theta,r)$ lies in the linear span of $(1,0)^T$ for all $V \in \cup_{k \geq 0} \mathcal{F}_k$, or $V(\theta,r)$ lies in the linear span of $(0,1)^T$ for all $V \in \cup_{k \geq 0} \mathcal{F}_k$. It follows that the strong bracket condition doesn't hold at any point $(\theta,r) \in M$.  
\eex 

In the previous example, the origin had to be excluded from $M$ in order to ensure that the unit circle is globally accessible. It could be interesting to determine whether there are PDMPs for which conditions \textbf{(i)} and \textbf{(ii)} are satisfied, the strong bracket condition nowhere holds, and $M$ is simply connected.  

As illustrated by the following example, the strong bracket condition does not imply condition \textbf{(i)}, not even if the vector fields are analytic.  

\bex \rm \label{ex:nonzero}
On the two-dimensional torus $\TT^2 = \RR^2 / \ZZ^2$, we switch between $F^0(x,y) = (1,0)^T$ and $F^1(x,y) = (0,1+\epsilon \sin(2 \pi x))^T$, where $\epsilon > 0$ is small. Any point in $\TT^2$ can then be reached from any starting point.  For $\alpha \in \RR$, we have 
$$ 
\alpha F^0(x,y) + (1-\alpha) F^1(x,y) = (\alpha, (1 -\alpha) (1+\epsilon \sin(2 \pi x)))^T, 
$$ 
which is never zero. However, 
$$ 
[F^0, F^1](x,y) = (0,-\epsilon 2 \pi \cos(2 \pi x))^T, 
$$ 
so the vectors $[F^0,F^1](0,0)$ and $F^0(0,0) - F^1(0,0) = (1,-1)^T$ span the tangent space at $(0,0)$, and the strong bracket condition is satisfied. 
\eex

\section{Applications}
\label{sec:app}
In this section, we give some applications of Theorem \ref{main} in the context of population models with an extinction set. For a general framework on Markov models with an extinction set, the reader is referred to \cite{B18}. Here we only give the results we will use in the specific context of PDMP on a compact set (see e.g \cite{BL16} or \cite{BS17}).
\subsection{Stochastic persistence}  
In this section, we assume that there exists a closed subset $M_0$ of $M$ which is invariant for the process : $X_t \in M_0$ if and only if $X_0 \in M_0$. The set $M_0$ will be referred to as the extinction set. We set $M_+ = M \setminus M_0$ and denote by $\mathcal{D}$ (resp. $\mathcal{D}^2$) the domain of the generator $\mathcal{L}$ defined in~\eqref{eq:defcL} (resp. the set of functions in the domain such that $f^2$ is also in  $\mathcal{D}$). We also let $\Gamma$ denote the carr{\'e} du champ operator on $\mathcal{D}^2$ : $\Gamma f = \mathcal{L}f^2 - 2 f \mathcal{L}f$, which  acts on functions $f \in \mathcal{D}^2$ as $$ \Gamma f = \sum_{j \in E} a_{ij}(x) \left(f^j(x) - f^i(x)\right)^2.$$

\begin{definition}
\label{def:pers}
We say that the process $Z$ is persistent if there exist continuous functions $V : M_+ \times E \to \RR_+$ and $H : M \times E \to \RR$ such that 
\begin{enumerate}
\item $\lim_{x \to M_0} V(x,i) = + \infty$,
\item For any compact set $K \subset M_+ \times E$, there exists $V_K \in \mathcal{D}^2$ such that $V|_K=V_K|_K$ and $({\cal L}V_K)|_K=H|_K$, 
\item There exists $\Delta > 0$ such that for all $t > 0$, $| V(Z_t) - V(Z_{t_-})| \leq \Delta$,
\item There exists $C > 0$ such that for any compact set $K \subset M_+$, $\| \Gamma(V_K)|_K\|_{\infty} \leq C$,
\item For any ergodic probability measure $\mu$ of $Z$ supported on $M_0 \times E$, one has $\mu H < 0$.
 \end{enumerate}
\end{definition}

The following theorem is an immediate consequence of \cite[Theorem 4.10]{B18} and Theorem \ref{main}.
\bthm
Assume that conditions \emph{\textbf{(i)}} and \emph{\textbf{(ii)}} hold, that $Z$ is persistent and that $e^{\star}$ is accessible from $M_+$. Then $Z$ admits a unique invariant probability measure $\Pi$ on $M_+ \times E$ and there exist $\theta, C, \gamma > 0$ such that  for all $t \geq 0$ and for all $(x,i) \in M_+ \times E$, $$ \| P_t((x,i),\cdot) - \Pi \|_{TV} \leq C \left( 1 + e^{\theta V(x,i)} \right) e^{- \gamma t}.$$
\label{thm:pers}
\ethm
\subsection{Lotka-Volterra in random environment}
\label{sec:LV}
In this section, we consider the competitive Lotka-Volterra model in a fluctuating environment studied in \cite{BL16} and show how our method can be used to improve one of their results.
 More precisely, for $i \in \{0,1\}$, let $F^i$ be defined as 
\begin{equation}\label{e:LVc}
F^i(x,y) =  \begin{pmatrix}
\alpha_i x (1 - a_i x - b_i y)\\
\beta_i y (1 - c_i x - d_i y)
\end{pmatrix},
\end{equation} 
with $\alpha_i, \beta_i, a_i, b_i, c_i, d_i > 0$. For $\eta > 0$ small enough, the flows $\varphi_t^i$ leave positively invariant the compact set $M = \{ (x,y) \in \RR_+^2 \: : \eta \leq x + y \leq 1/ \eta \}$, and the extinction set $M_0$ is the union of $M_0^1 = \{(x,y) \in M \: : x = 0 \}$ and  $M_0^2 = \{(x,y) \in M \: : y = 0 \}$.    It is shown in \cite{BL16} that the long-term behavior of the process $(Z_t)_{ t \geq 0}=(X_t,Y_t,I_t)_{ t \geq 0}$ is determined by the sign of the invasion rates : $$\Lambda_y = \int \beta_i (1 - c_i x) \mathrm{d} \mu(x,i),$$ and $$ \Lambda_x  = \int \alpha_i  (1 - b_i y) \mathrm{d} \hat{\mu} (y,i),$$ where $\mu$ and $\hat{\mu}$ are the unique invariant probability measures of the process $Z$ restricted to $M_0^2$ and $M_0^1$, respectively. It is not hard to construct functions $V : M_+ \times E \to \RR$ and $H : M \times E \to \RR_+$ satisfying assumptions $1.$ to $4.$ of Definition \ref{def:pers}, such that $V(x,y,i)$ co\"incides with $-\log(x)$ in a neighborhood of $M_0^1$ and with $-\log(y)$ in a neighborhood of $M_0^2$, and such that $H(x,y,i)$ co\"incides with $ \alpha_i  (1 -a_i x - b_i y )$ in a neighbourhood of $M_0^1$ and with $\beta_i (1 - c_i x - d_i y)$ in a neighborhood of $M_0^2$ (see e.g \cite[Section 5]{B18} or \cite[Section 5]{BS17}). Then, one can check that  $\Lambda_x= - \hat{\mu} H$ and $\Lambda_y = - \mu H$, so that $Z$ is persistent if and only if $\Lambda_x > 0 $ and $\Lambda_y >0$.

It is shown in \cite{BL16} that if $\Lambda_x > 0 $ and $\Lambda_y >0$, then the process admits a unique invariant probability measure $\Pi$ in $M_+ \times E$. But to show the convergence in total variation of the law of $Z_t$ toward $\Pi$, the authors needed to check that the strong bracket condition is satisfied at some accessible point. They proved, except in the particular case where $\frac{\beta_0 \alpha_1}{\alpha_0 \beta_1}= \frac{a_0 c_1}{c_0 a_1}= \frac{b_0 d_1}{d_0 b_1}$, that this condition holds by using a formal calculus program.  Thanks to Theorem \ref{thm:pers}, we withdraw this condition, and give an easier proof for the convergence in total variation.

In \cite{BL16}, of particular importance is the study of the averaged vector fields $F^s := sF^1 + (1-s) F^0$, for $s \in [0,1]$. The vector field $F^s$ is still a competitive Lotka-Volterra system of the form \eqref{e:LVc}, with coefficients   $\alpha_s, \beta_s, a_s, b_s, c_s, d_s$ that are barycentric combinations of the coefficients appearing in $F^0$ and $F^1$. The dynamics of the deterministic system generated by $F^s$ depends on the position of $s$ with respect to the two following (possibly empty) intervals: $$ I = \{ s \in (0,1) \: : a_s > c_s \}$$ and $$ J = \{ s \in (0,1) \: : b_s > d_s \}.$$ There are four regions of interest :
\begin{itemize}
\item $s \in (\overline{I})^c \cap (\overline{J})^c$ : the equilibrium $(1/a_s,0)$ is a global attractor for solutions with $x_0 \neq 0$;
\item $s \in I \cap J$ : the equilibrium $(0,1/b_s)$ is a global attractor for solutions with $y_0 \neq 0$;
\item $s \in I \cap (\overline{J})^c$ : $F^s$ admits a unique G.A.S. equilibrium $e_s \in M_+$;
\item $s \in (\overline{I})^c \cap J$ : $F^s$ admits a unique equilibrium $e_s \in M_+$, which is a saddle whose stable manifold separates the basins of attraction of $(1/a_s,0)$ and $(0,1/b_s)$.
\end{itemize}
Here, $(\overline{I})^c$ and $(\overline{J})^c$  stand for the complement of the closure of $I$ and $J$, respectively.  The following proposition is a consequence of \cite[Proposition 2.3 and Theorem 4.1]{BL16}.

\bprop
Assume $\Lambda_y > 0$. Then $I \neq \emptyset$ and there exists a point $m$ accessible from $M_+$ such that the weak bracket condition holds at $m$.
\label{prop:lambda>0}
\eprop

From this proposition, we can derive the next lemma:

\blem
Assume $\Lambda_y > 0$. Then there exists $s \in [0,1]$ such that $F^s$ admits an equilibrium $e_s \in M_+$ which is accessible from $M_+$. In particular, condition \emph{\textbf{(i)}} holds.
\label{lem}
\elem

This lemma combined with Proposition \ref{prop:lambda>0} and Theorem \ref{thm:pers} implies the following corollary, which slightly improve \cite[Theorem 4.1 - (iv)]{BL16}

\bcor
Assume $\Lambda_y > 0$ and $\Lambda_x >0$. Then there exist $C, \gamma, \theta > 0$ such that for all $t \geq 0$ and for all $(x,y,i) \in M_+ \times E$, $$ \| P_t((x,i),\cdot) - \pi \|_{TV} \leq C \left( 1 + \frac{1}{\|x\|^{\theta}} + \frac{1}{\|y\|^{\theta}} \right) e^{- \gamma t}.$$ 
\ecor
\bigskip

\prf \textbf{of Lemma \ref{lem}}
Since $\Lambda_y > 0$, $I$ is nonempty by Proposition \ref{prop:lambda>0}. Then we have three cases: either $I \cap J^c$ is nonempty,  or $I$ is a strict subset of $J$  or $I = J$. We prove the lemma in these three cases. Assume first that $I \cap J^c \neq \emptyset$ and take $s \in I \cap J^c$. Then $F^s$ admits a G.A.S. equilibrium $e_s \in M_+$, in particular it is accessible. Assume now that $I$ is a strict subset of $J$. In particular, $I^c \cap J$ and $I \cap J$ are nonempty. Pick $s \in I^c \cap J$, then  $F^s$ admits a unique equilibrium $e_s \in M_+$, which is a saddle whose stable manifold $W_s$ separates the basins of attraction of $(1/a_s,0)$ and $(0,1/b_s)$. We show that $e_s$ is accessible. Choose a point $(x,y) \in M_+$. Then, if $(x,y)$ is above $W_s$, follow the flow $\varphi^0$. As the resulting trajectory converges to $(1/a_0,0)$, it needs to cross $W_s$. If $(x,y)$ is below $W_s$, one can find a trajectory leading to $(0,1/b_u)$ for some $u \in I \cap J$. In particular, this trajectory also crosses $W_s$. As $e_s$ is also accessible from every point in $W_s$, it is accessible from everywhere in $M_+$. Finally, assume that $I=J=(s_1,s_2)$. Then the vector field $F^{s_1}$ is of the form 
$$F^{s_1}(x,y) =  \begin{pmatrix}
\alpha x (1 - a x - b y)\\
\beta y (1 - a x - b y)
\end{pmatrix},$$ with $a = a_{s_1} = c_{s_1}$ and $b= b_{s_1}=d_{s_1}$. In particular, the line $y = 1/b(1-ax)$ is composed of equilibria of $F^{s_1}$. Moreover, $(1/a_0,0)$ and $(1/a_1,0)$ lie on opposite sides of this line. Now we know by Proposition \ref{prop:lambda>0} that there exists an accessible point $m \in M_+$. Hence, depending on the position of $m$ with respect to the line $y = 1/b(1-ax)$, follow either $\varphi^0$ or $\varphi^1$ in order to cross the line when starting at $m$. Then the point where the line is crossed is accessible from $m$ and therefore from $M_+$.\qed

\subsection{Epidemiological models : SIS in dimension 2}

In this section we discuss an application of Theorem \ref{thm:pers} to an SIS model with two groups and two environments, as studied in \cite[Section 4]{BS17}. We look at random switching between differential equations on $[0,1]^2$ having the form
\beq
\label{eq:LY}
\frac{dx_i}{dt} = (1-x_i) (\sum_{j = 1}^d C^k_{ij} x_j) - D^k_i x_i\, ,   \: i = 1, 2,
\eeq
where for $k \in E =\{0,1\}$, $C^k = (C^k_{ij})$ is an irreducible matrix with nonnegative entries and $D^k_i > 0$. Let $A^k = C^k - \mathsf{diag}(D^k)$ and let  $\lambda(A^k)$ denote the largest real part of the eigenvalues of $A^k$. Then, we have the following result due to Lajmanovich and Yorke.
\bthm[Lajmanovich and Yorke, \cite{LajYorke}]
 \label{th:LajYorke} 
 If $\lambda(A^k) \leq 0,$  $0$ is a G.A.S equilibrium for the semiflow induced by (\ref{eq:LY}) on $[0,1]^2.$
If $\lambda(A^k) > 0,$ there exists another equilibrium $x^*_k \in (0,1)^2$ whose basin of attraction is $[0,1]^2 \setminus \{0\}.$
\ethm

\blem
\label{lem:LY}
Assume that 
\begin{enumerate}
\item  $\lambda(A^0)<0$ and $\lambda(A^1)<0$,
\item There exists $s \in (0,1)$ such that $\lambda(A^s)>0$, where $A^s = s A^1 + (1-s) A^0$.
\end{enumerate}
Then conditions \emph{\textbf{(i)}} and \emph{\textbf{(ii)}} are satisfied.
\elem
An example where the assumptions of this lemma hold can be found in \cite[Example 4.7]{BS17}. If the assumptions of Lemma \ref{lem:LY} hold, Corollary 2.14 and Section 5 in \cite{BS17} imply that $Z$ is persistent provided the switching occurs sufficiently often. In that case, we get by Theorem \ref{thm:pers} the convergence in total variation to a unique invariant probability measure.  Compare this to \cite[Theorem 4.11]{BS17}, which only gives convergence in a certain Wasserstein distance. Note that the conclusion of Lemma \ref{lem:LY} is no longer true in general if $\lambda(A^0) > 0$ and  $\lambda(A^1) > 0$. An easy counterexample is when the two equilibria $x^*_0$, $x^*_1$ given by Theorem \ref{th:LajYorke} co\"incide (see e.g. \cite[Example 4.10]{BS17}). In that case, condition \textbf{(i)} is satisfied but condition \textbf{(ii)} obviously is not.

\prf \textbf{of Lemma \ref{lem:LY}}
  For $k \in E$, we let $F^k$ denote the vector field given by the right hand side of \eqref{eq:LY}. It is readily seen that for $s \in (0,1)$, the vector field $F^s = s F^1 + (1-s)F^0$ is of the same form as $F^0$ and $F^1$, with matrix $C^s=s C^1 + (1-s) C^0$ and vector $D^s = s D^1 + (1-s) D^0$. As a consequence, since there exists $s \in (0,1)$ such that $\lambda(A^s)>0$, Theorem \ref{th:LajYorke} implies that condition \textbf{(i)}  is satisfied at some point $x^*_s \in (0,1)^2$, and we even have $F^s(x^*_s)=0$. Moreover, since  $\lambda(A^0)<0$ and $\lambda(A^1)<0$, the first part of Theorem \ref{th:LajYorke} implies that neither $F^0$ nor $F^1$ can vanish at $x_s^*$. In particular, $F^0(x_s^*)$ and $F^1(x_s^*)$ are collinear and of opposite direction. For $k \in \{0,1\}$ let $\gamma^k(x_s^*)$ denote the positive orbit of $x_s^*$ under $F^k$. Due to the first part of Theorem \ref{th:LajYorke}, $\gamma^0(x_s^*)$  is a curve linking $x_s^*$ and $0$. To obtain a contradiction, assume that condition \textbf{(ii)} is not satisfied. Then $F^0$ and $F^1$ are collinear and of opposite direction on $\gamma^0(x_s^*)$. We have for all $x \in  \gamma^0(x_s^*)$ that $x_s^* \in  \gamma^1(x)$, meaning that for all $\eps >0$, one can find $x$ with $\|x\| < \eps$ and $t>0$ such that $\| \varphi_t^1(x) \| = \| x_s^* \|$. This is in contradiction with the fact that $0$ is a G.A.S equilibrium for $F^1$, hence condition \textbf{(ii)} holds as well.
\qed

\section{Proof of Proposition~\ref{prop:submersion_at_E}} \label{sec:proof}
To prove Proposition~\ref{prop:submersion_at_E}, we will use \cite[Theorem 4.1]{BMZIHP} that we quote here.

\bthm[Bena\"im -- Le Borgne -- Malrieu -- Zitt]
Let $x$ be a point of $M$ at which the weak bracket condition holds. Then, there exists $m \geq d$, $\bi=(i_1, \ldots, i_m) \in E^m$ and $\bu = (u_1,\ldots,u_m) \in \RR_+^m$ such that the map $ \bv \to \bphi_{\bv}^{\bi}(x)$ is a submersion at $\bu$.
\label{thmwbsub}
\ethm
 
The following proposition is the key point of the proof :

\bprop
Under the hypothesis of Theorem \ref{main}, there exist $s > 0$, $i \in E$, $\bi=(i_1, \ldots, i_n) \in E^n$ and $\bu = (u_1,\ldots,u_n) \in \RR_+^n$ with $s > u_1 + \ldots + u_n $ such that the map $ \Psi : D^s_{n+1} \to \RR^d  $, $ (\bv,t) \to \varphi_{s-\sum v_i-t}^{i} \circ \bphi_{\bv}^{\bi}(e^{\star})$ is a submersion at $(\bu,0)$.
\label{mainp}
\eprop

This proposition remains valid if we replace $e^{\star}$ by any point in $M$ from which one can access a point $x^*$ where the weak bracket condition holds.  In particular, it is independent of our assumption that $e^{\star}$ is an equilibrium of a vector field of the form $\sum \alpha_i F^i$. The proposition is a consequence of the two lemmas we give now.
\blem
Suppose that there exists a point $x^*$ accessible from $e^{\star}$ such that the weak bracket condition holds at $x^*$. Then there exists $(\bar{\bi},\bar{\bu})$ such that the weak bracket condition holds at $\bphi_{\bar{\bu}}^{\bar{\bi}}(e^{\star})$.
\elem 

\prf
By Proposition \ref{acces},  $x^*$ is accessible from $e^{\star}$ if and only if $x^*\in \overline{\gamma^+(e^{\star})}$. By continuity of the determinant and regularity of the vector fields, the weak bracket condition is an open condition. Thus if it holds at a point of $\overline{\gamma^+(e^{\star})}$, it also holds at a point in $\gamma^+(e^{\star})$, hence the result.
\qed
Thanks to this lemma, we assume from now on that  there exist $\bar{\bi}=(\bar{i}_1,\ldots,\bar{i}_p)$ and $\bar{\bu}=(\bar{u}_1,\ldots,\bar{u}_p)$ such that $x^*=\bphi_{\bar{\bu}}^{\bar{\bi}}(e^{\star})$. Since $x^*$ satisfies the weak bracket condition,  Theorem \ref{thmwbsub} implies that there exists $m \geq d$, $\bi=(i_1, \ldots, i_m) \in E^m$ and $\bu = (u_1,\ldots,u_m) \in \RR_+^m$ such that the map $ \psi : \bv \to \bphi_{\bv}^{\bi}(x^*)$ is a submersion at $\bu$. We denote $\bi_- = (i_1, \ldots, i_{m-1})$ and $\bv_- =(v_1, \ldots, v_{m-1})$, and for all $s > 0$, we define the map $\Psi^s : D^s_{m+p} \to \RR^d$ by 

$$\Psi^s : (\bv_-,\bar{\bv},t) \to \varphi_{s-(v_1+\ldots+v_{m-1}+\bar{v}_1+\ldots+\bar{v}_p +t)}^{i_m} \circ \Phi_{\bv_-}^{\bi_-} \circ \bphi_{\bar{\bv}}^{\bar{\bi}}(e^{\star}).$$
 We also let $\sigma_{(\bv_-,\bar{\bv})}^t = v_1+\ldots+v_{m-1}+ \bar{v}_1+ \ldots + \bar{v}_p+t$.  Note that in particular, $$\Psi^s (\bv_-,\bar{\bu},t)= \varphi_{s-\sigma_{(\bv_-,\bar{\bu})}^t}^{i_m} \circ \Phi_{\bv_-}^{\bi_-} (x^*)=\psi(\bv_-,s-\sigma_{(\bv_-,\bar{\bu})}^t)$$ for all $(\bv_-, \bar{\bu}, t) \in D^s$.  With this property, the next lemma is straightforward :

\blem
For all $k \in \{1, \ldots, m-1\}$, for all $(\bv_-,\bar{\bu},t) \in D^s_{m+p}$, one has

$$\frac{\partial \Psi^s}{\partial v_k}(\bv_-,\bar{\bu},t) = - \frac{\partial \psi}{\partial v_m}(\bv_-,s-\sigma_{(\bv_-,\bar{\bu})}^t) +   \frac{\partial \psi}{\partial v_k}(\bv_-,s-\sigma_{(\bv_-,\bar{\bu})}^t),
$$ and
$$
\frac{\partial \Psi^s}{\partial t}(\bv_-,\bar{\bu},t)= - \frac{\partial \psi}{\partial v_m}(\bv_-,s-\sigma_{(\bv_-,\bar{\bu})}^t).
$$ In particular, setting $s = u_1 + \ldots + u_m + \bar{u}_1 + \ldots + \bar{u}_p$ and $t=0$, one gets 
\beq
\frac{\partial \Psi^s}{\partial v_k}(\bu_-,\bar{\bu},0) = - \frac{\partial \psi}{\partial v_m}(\bu) +   \frac{\partial \psi}{\partial v_k}(\bu),
\label{partialk}
\eeq and 
\beq
\frac{\partial \Psi^s}{\partial t}(\bu_-,\bar{\bu},0)= - \frac{\partial \psi}{\partial v_m}(\bu).
\label{partialt}
\eeq
\elem
 
%
%
%
%
\bigskip

\prf \textbf{of Proposition \ref{mainp}}

For $s=u_1+\ldots+u_m+\bar{u}_1+\ldots+\bar{u}_p$, equalities \eqref{partialk} and \eqref{partialt} proves that the rank of the family of vectors $(\frac{\partial \Psi^s}{\partial v_1}(\bu_-,\bar{\bu},0), \ldots, \frac{\partial \Psi^s}{\partial v_{m-1}}(\bu_-,\bar{\bu},0),\frac{\partial \Psi^s}{\partial t}(\bu_-,\bar{\bu},0))$ is the same as the family of vectors $(\frac{\partial \psi}{\partial v_k}(\bu), 1 \leq k \leq m)$. But since $\psi$ is a submersion at $\bu$, this rank is $d$, showing that $\Psi^s$ is also a submersion at point $(\bu_-,\bar{\bu},0)$.\qed

We can now pass to the main part of the proof of Proposition~\ref{prop:submersion_at_E}.

\bigskip

\prf \textbf{of Proposition~\ref{prop:submersion_at_E}}

We first construct a function $\bar{\Psi}$ and then verify that it is indeed a submersion.  By Proposition \ref{mainp}, there exist $s > 0$, $\bi=(i_1, \ldots, i_n,i_{n+1}) \in E^{n+1}$ and $\bu = (u_1,\ldots,u_n) \in \RR_+^n$ such that the map $ \Psi : (\bv,t) \to \varphi_{s-\sum v_i-t}^{i_{n+1}} \circ \bphi_{\bv}^{\bi}(e^{\star})$ is a submersion at $(\bu,0)$. In the sequel, we denote by $\Psi  (\bv,t)$ the map given by $\Psi  (\bv,t)(x) = \varphi_{s-\sum v_i-t}^{i_{n+1}} \circ \bphi_{\bv}^{\bi}(x)$. We define the map $\overline{\Psi}$ on $D^s_{n+N}$ with values in $\RR^d$ by $$ \overline{\Psi}(\bv, \bar{\bv}) \to \varphi_{s-\sum v_i-\sum \bar{v_i}}^{i_{n+1}} \circ \bphi_{\bv}^{\bi} \circ \bphi_{\bar{\bv}}^{\bar{\bi}}( e^{\star}),$$ where $\bar{\bi} = (1,2,\ldots, N)$. Then with the previous notation, $\overline{\Psi}(\bv, \bar{\bv}) = \Psi( \bv, \sum \bar{v}_i) \circ  \bphi_{\bar{\bv}}^{\bar{\bi}}( e^{\star})$.  Now, we show that the map $\overline{\Psi}$ is a submersion at $(\bu,0)$ --- here, $0$ denotes the zero vector in $\RR^N$. For all $1 \leq k \leq n$, 
\beq
\frac{\partial \overline{\Psi}}{\partial v_k}(\bv,\bar{\bv})=\frac{\partial \Psi}{\partial v_k}(\bv, \sum \bar{v}_i) \circ  \bphi_{\bar{\bv}}^{\bar{\bi}}( e^{\star}),
\label{firstpartial}
\eeq
and for all $1 \leq k \leq N$,
\beq
\frac{\partial \overline{\Psi}}{\partial \bar{v}_k}(\bv,\bar{\bv})=\frac{\partial \Psi}{\partial t}(\bv, \sum \bar{v}_i) \circ  \bphi_{\bar{\bv}}^{\bar{\bi}}( e^{\star}) + D\Psi( \bv, \sum \bar{v}_i)(\bphi_{\bar{\bv}}^{\bar{\bi}}(e^{\star}) ) \frac{\partial \bphi_{\bar{\bv}}^{\bar{\bi}}}{\partial \bar{v}_k}(e^{\star}). 
\eeq
Now, since each $\varphi^i_v$ is the identity at $v = 0$ and $\partial_v \varphi^i_v(x) = F^i(\varphi^i_v(x))$, one can easily show that 
\beq
\left. \frac{\partial \bphi_{\bar{\bv}}^{\bar{\bi}}}{\partial \bar{v}_k}(e^{\star}) \right|_{\bar{\bv}=0} = F^k(e^{\star}).
\eeq
In particular, since $\bphi_{\bar{\bv}}^{\bar{\bi}}( e^{\star})=e^{\star}$ when $\bar{\bv}=0$, $$\frac{\partial \overline{\Psi}}{\partial \bar{v}_k}(\bv,0)=\frac{\partial \Psi}{\partial t}(\bv, 0)( e^{\star})+ D\Psi( \bv, \sum \bar{v}_i)(e^{\star} ) F^k(e^{\star}),$$
which, due to condition \textbf{(i)} implies that 
\beq
 \sum_{k=1}^N \alpha_k \frac{\partial \overline{\Psi}}{\partial \bar{v}_k}(\bv,0)=\frac{\partial \Psi}{\partial t}(\bv, 0)( e^{\star}).
\label{average}
\eeq
Thus, \eqref{firstpartial} and \eqref{average} evaluated at $\bv =\bu$ and $\bar{\bv}=0$ yield $$\mathrm{rank} \left( \frac{\partial \overline{\Psi}}{\partial v_k}(\bu,0),\frac{\partial \overline{\Psi}}{\partial \bar{v}_k}(\bu,0) \right) \geq \mathrm{rank} \left(\frac{\partial \Psi}{\partial v_k}(\bu, 0),\frac{\partial \Psi}{\partial t}(\bu, 0) \right)=d, $$
where the last equality is due to Proposition \ref{mainp}. This finishes the proof. \qed

\section*{Acknowledgments}
This work was supported by the SNF grants $200021-175728$ and $200021L\_169691$. We thank an anonymous referee for useful remarks. 

%
%
 
\bibliographystyle{amsalpha}
\bibliography{RandomSwitch}

\providecommand{\bysame}{\leavevmode\hbox to3em{\hrulefill}\thinspace}
\providecommand{\MR}{\relax\ifhmode\unskip\space\fi MR }
\providecommand{\MRhref}[2]{%
  \href{http://www.ams.org/mathscinet-getitem?mr=#1}{#2}
}
\providecommand{\href}[2]{#2}
\begin{thebibliography}{BLBMZ15}

\bibitem[BCL17]{BCL17}
M.~Bena{\"\i}m, F.~Colonius, and R.~Lettau, \emph{Supports of invariant
  measures for piecewise deterministic {M}arkov processes}, Nonlinearity
  \textbf{30} (2017), no.~9, 3400--3418. \MR{3694261}

\bibitem[Ben18]{B18}
M.~Bena{\"\i}m, \emph{Stochastic persistence},
  \url{https://arxiv.org/abs/1806.08450}.

\bibitem[BH12]{BH12}
Y.~Bakhtin and T.~Hurth, \emph{Invariant densities for dynamical systems with
  random switching}, Nonlinearity \textbf{25} (2012), no.~10, 2937--2952.
  \MR{2979976}

\bibitem[BHLM18]{torus}
Yuri Bakhtin, Tobias Hurth, Sean~D Lawley, and Jonathan~C Mattingly,
  \emph{Smooth invariant densities for random switching on the torus},
  Nonlinearity \textbf{31} (2018), no.~4, 1331.

\bibitem[BL16]{BL16}
M.~Bena{\"\i}m and C.~Lobry, \emph{Lotka {V}olterra in fluctuating environment
  or ``how switching between beneficial environments can make survival
  harder''}, Ann. Appl. Probab. \textbf{26} (2016), no.~6, 3754--3785.

\bibitem[BLBMZ12]{Le_Borgne}
Michel Bena{\"{\i}}m, St{\'e}phane Le~Borgne, Florent Malrieu, and
  Pierre-Andr{\'e} Zitt, \emph{Quantitative ergodicity for some switched
  dynamical systems}, Electron. Commun. Probab. \textbf{17} (2012), no. 56, 14.
  \MR{3005729}

\bibitem[BLBMZ15]{BMZIHP}
M.~Bena{\"\i}m, S.~Le~Borgne, F.~Malrieu, and P.~A. Zitt, \emph{Qualitative
  properties of certain piecewise deterministic {M}arkov processes}, Ann. Inst.
  Henri Poincar{\'e} Probab. Stat. \textbf{51} (2015), no.~3, 1040--1075.
  \MR{3365972}

\bibitem[BS17]{BS17}
M.~Bena{\"\i}m and E.~Strickler, \emph{Random switching between vector fields
  having a common zero}, \url{https://arxiv.org/abs/1702.03089}.

\bibitem[CH15]{Cloez}
B.~Cloez and M.~Hairer, \emph{Exponential ergodicity for {M}arkov processes
  with random switching}, Bernoulli \textbf{21} (2015), 505--536.

\bibitem[Dav84]{Dav84}
M.~H.~A. Davis, \emph{Piecewise-deterministic {M}arkov processes: a general
  class of nondiffusion stochastic models}, J. Roy. Statist. Soc. Ser. B
  \textbf{46} (1984), no.~3, 353--388. \MR{790622}

\bibitem[LLC17]{li17}
Dan Li, Shengqiang Liu, and Jing'an Cui, \emph{Threshold dynamics and
  ergodicity of an {SIRS} epidemic model with markovian switching}, Journal of
  Differential Equations \textbf{263} (2017), no.~12, 8873--8915.

\bibitem[Nua06]{Nualart}
David Nualart, \emph{The {M}alliavin calculus and related topics}, second ed.,
  Probability and its Applications (New York), Springer-Verlag, Berlin, 2006.
  \MR{2200233 (2006j:60004)}

\bibitem[SJ72]{SJ72}
H{\'e}ctor~J. Sussmann and Velimir Jurdjevic, \emph{Controllability of
  nonlinear systems}, J. Differential Equations \textbf{12} (1972), 95--116.

\end{thebibliography}
\end{document}